\documentclass[preprint,12pt]{elsarticle}

\usepackage{amssymb}
\usepackage{amsthm}
\usepackage{amsmath}

\journal{Journal of Number Theory}
\newtheorem{theorem}{Theorem}
\newtheorem{lemma}{Lemma}
\newtheorem*{corollary}{Corollary}
\numberwithin{equation}{section}
\newproof{pot1}{Proof of Theorem \ref{theorem1}}
\newproof{pot2}{Proof of Theorem \ref{theorem2}}
\begin{document}

\begin{frontmatter}



\title{Asymptotics of class number and genus for abelian extensions of an algebraic function field}


\author{Kenneth Ward\corref{cor1}}
\cortext[cor1]{Corresponding author}
\address{Department of Mathematics, Oklahoma State University, Stillwater, OK 74078 USA}
\ead{kward@math.okstate.edu}

\begin{abstract}
Among abelian extensions of a congruence function field, an asymptotic relation of class number and genus is established: namely, for such extensions with class number $h$, genus $g$, and field of constants $\mathbb{F}$, that $\ln h \sim g \ln |\mathbb{F}|$.  The proof is completely classical, employing well known results from congruence function field theory. This gives an answer to a question of E. Inaba.
\end{abstract}


\end{frontmatter}


\section{Introduction}
\label{introduction}

Let $K$ be a congruence function field with genus $g_K$ and class number $h_K$. The study of the asymptotic behavior of class number and genus for congruence function fields dates to a result of E. Inaba \cite{In}, which established, for a natural number $m$, that among congruence function fields $K$ with a fixed choice of finite constant field $\mathbb{F}_q$ and an element $x \in K$ that satisfies $[K:\mathbb{F}_q(x)]\leq m$, \begin{equation}\label{eq1}\lim_{g_K \rightarrow \infty} \frac{\ln h_K}{g_K \ln |\mathbb{F}_q|} = 1.\end{equation} In his paper Inaba remarked that he was not aware of whether this relation remains true if the bound involving $m$ is removed. As noted by K. Iwasawa in Inaba's article, the requirement that $m$ be fixed resembles R. Brauer's first result on the Brauer-Siegel theorem for algebraic number fields \cite{Br1}. Results similar to that of Inaba also appear in the work of I. Luthar and S. Gogia \cite{GoLu} and M. Tsfasman \cite{Ts}.

Much later, M. Madan and D. Madden \cite{MaMa} noted, for congruence function fields $K$ with a fixed choice of constant field $\mathbb{F}_q$ and an element $x\in K\backslash \mathbb{F}_q$, that Inaba's method yields \begin{equation} \label{eq2} \lim_{\frac{[K:\mathbb{F}_q(x)]}{g_K}  \rightarrow 0}  \frac{\ln h_K}{g_K \ln |\mathbb{F}_q|} = 1. \end{equation} One may observe that (\ref{eq2}) loosely resembles the condition required in R. Brauer's second paper \cite{Br2} on the Brauer-Siegel theorem, with an exception: in addition to requiring for an extension $L$ of the rational numbers with discriminant $d$ that $[L:\mathbb{Q}] / \ln d$ tends to zero, it was necessary for Brauer to assume that the extension $L / \mathbb{Q}$ be normal. It is difficult to surmount both of these requirements in the case of number fields as a result of the connection to the Riemann hypothesis \cite{La2}. However, Brauer's result may be extended to abelian number fields without any relative requirement on discriminant growth; for example, see the concise argument of S.R. Louboutin \cite{Lou}.

The objective of this paper is to establish the analogue to Brauer's theorem for finite abelian extensions of \emph{any} choice of congruence function field. In fact, P. Lam-Estrada and G.D. Villa-Salvador \cite{LaVi} have already noted that, by a result of D. Hayes \cite{Ha}, the relation (\ref{eq1}) holds among the cyclotomic extensions of the field of rational functions $\mathbb{F}_q(T)$. The objective is met by means of two theorems. For what follows, let $\mathbb{F}_K$ denote the constant field of $K$.

\begin{theorem}\label{theorem1} Let $K$ be a congruence function field. It holds that \begin{equation}\notag \liminf_{g_K \rightarrow \infty} \frac{\ln h_K}{g_K \ln \left|\mathbb{F}_K\right|} \geq 1. \end{equation} \end{theorem}

The bound attained in the proof of Theorem 1 is effective. Furthermore, Theorem 1 makes no requirement that $K$ be a finite abelian extension of a congruence function field.  It is for the proof of the upper bound that abelian structure is essential. 

\begin{theorem}\label{theorem2} Let $F$ be a congruence function field. Let $K$ be a finite abelian extension of $F$. It holds that \begin{equation}\notag \limsup_{g_K \rightarrow \infty} \frac{\ln h_K}{g_K \ln \left|\mathbb{F}_K\right|} \leq 1. \end{equation} \end{theorem} 

The method of proof for Theorem 2 in this paper is unique to the abelian case; indeed, the required properties are violated within the simplest class of non-abelian extensions of a congruence function field for which the genus may become large: finite, geometric, tamely ramified, and solvable extensions of $\mathbb{F}_q(T)$. This is a consequence of the possibility of slow growth of the genus \cite{GaStTh}. Also, unlike Theorem 1, the bound attained in the proof of Theorem 2 is ineffective.  As a corollary of Theorems 1 and 2, one obtains the main result of this paper. 

\begin{corollary} Let $F$ be a congruence function field. Let $K$ be a finite abelian extension of $F$. It holds that \begin{equation}\notag \lim_{g_K \rightarrow \infty} \frac{\ln h_K}{g_K \ln |\mathbb{F}_K|}=1.\end{equation}  \end{corollary}

\section{The lower bound}
\label{section2}

The proof of Theorem 1 proceeds as follows.

 \begin{enumerate} \item Count the number of monic irreducible polynomials of a given degree with coefficients in $\mathbb{F}_K$ via M\"{o}bius inversion; \item For $x \in K \backslash \mathbb{F}_K$, compare the number of places of a given degree in $K$ to the number of places of the same degree in $\mathbb{F}_K(x)$ via M\"{o}bius inversion and Riemann's hypothesis; \item Obtain a lower bound for the number of integral divisors of degree $2g_K$ in $K$ via the Riemann-Roch theorem. \end{enumerate}
 
This proof follows closely Inaba's original method in \cite{In}. The first step is a basic result in field theory \cite{La1}.

\begin{lemma}\label{lemma1} Let $x \in K \backslash \mathbb{F}_K$. For each $m \in \mathbb{N}$, let $\psi(m)$ be the number of monic irreducible elements of $\mathbb{F}_K[x]$ of degree in $x$ equal to $m$. Let $\mu$ be the M\"{o}bius function. It holds for each $m \in \mathbb{N}$ that \begin{equation}\notag \psi(m) = \frac{1}{m} \sum_{d|m} \mu\left(\frac{m}{d}\right)\left|\mathbb{F}_K\right|^d.\end{equation} \end{lemma}

The second step of the proof follows a method known to H. Reichardt \cite{Re}. The basic principle is as follows. For $K$, let $\mathbb{P}_K$ denote the collection of places and $d_K$ the degree function on divisors. For each $m \in \mathbb{N}$, let \begin{equation}\notag N_m = |\{\mathfrak{P} \in \mathbb{P}_K \;|\; d_K(\mathfrak{P})= m\}|. \end{equation} Also, let $s \in \mathbb{C}$ with $\text{Re}(s) > 1$  and $u = |\mathbb{F}_K|^{-s}$. One may write the zeta function $\zeta_K(s)$ of $K$ as \begin{equation} \label{euler} \zeta_K(s) =  \prod_{\mathfrak{P} \in \mathbb{P}_K} \left(1 - \frac{1}{|\mathbb{F}_K|^{d_K(\mathfrak{P})s}}\right)^{-1} = \prod_{k=1}^\infty \left(1 - u^k\right)^{-N_k}. \end{equation} Let $x \in K \backslash \mathbb{F}_K$. For $\mathbb{F}_K(x)$, let  $\mathbb{P}_0$ denote the collection of places, $d_0$ the degree function on divisors, $\zeta_0(s)$ the zeta function, and \begin{equation} \notag n_m = |\{\mathfrak{p} \in \mathbb{P}_0 \;|\; d_0(\mathfrak{p}) = m\}|. \end{equation} Application of the logarithmic derivative to both (\ref{euler}) and the analogous identity for $\zeta_0(s)$ yields that \begin{equation}\label{logarithmicderivative} \frac{\zeta'_K(s)}{\zeta_K(s)} - \frac{\zeta'_0(s)}{\zeta_0(s)}= -\ln |\mathbb{F}_K| \sum_{m=1}^\infty \left(\sum_{d|m} d \left(N_d-n_d\right) \right) u^m. \end{equation} Let $P_K(s) = (1 - u)(1 - |\mathbb{F}_K|u)\zeta_K(s)$. It is well known \cite{Deu} that there exist $\omega_1,...,\omega_{2g_K} \in \mathbb{C}$ with \begin{equation}\label{zetanumerator} P_K(s) = \prod_{i=1}^{2g_K} \left( 1 - \omega_i u \right). \end{equation} Furthermore, as $\text{Re}(s)>1$, one has that $P_K(s)=\zeta_K(s)/\zeta_0(s)$.  From (\ref{logarithmicderivative}) and (\ref{zetanumerator}), it follows that \begin{equation}\label{inversion} \sum_{d|m} d(N_d - n_d) = - \sum_{i=1}^{2g_K} \omega_i^m. \end{equation} By Riemann's hypothesis, it holds for each $i=1,...,2g_K$ that $|\omega_i| = |\mathbb{F}_K|^{\frac{1}{2}}$. By M\"{o}bius inversion, one then obtains from (\ref{inversion}) the following lemma.

\begin{lemma}\label{lemma2} Let $x \in K \backslash \mathbb{F}_K$. For each $m \in \mathbb{N}$, it holds that \begin{equation} \notag |N_m - n_m| \leq 4g_K |\mathbb{F}_K|^{\frac{m}{2}}. \end{equation} \end{lemma}

\begin{pot1} For a divisor class $C$ of $K$, let $l_K(C)$ denote the dimension over $\mathbb{F}_K$ of the Riemann-Roch space for any element of $C$. If $C$ is of degree equal to $2g_K$, the Riemann-Roch theorem gives that $l_K(C)=g_K+1$. Thus the number of integral divisors $A_{2g_K}$ of $K$ of degree $2g_K$ satisfies \begin{equation} \label{divisor2gcount} A_{2g_K} = h_K \left( \frac{|\mathbb{F}_K|^{g_K+1}-1}{|\mathbb{F}_K|-1} \right). \end{equation} By (\ref{divisor2gcount}) and Lemmas \ref{lemma1} and \ref{lemma2}, one obtains that \begin{align}\label{classnumberlowerbound1} \notag h_K \left( \frac{|\mathbb{F}_K|^{g_K+1}-1}{|\mathbb{F}_K|-1} \right)& \geq N_{2g_K} \geq n_{2g_K} - 4g_K |\mathbb{F}_K|^{g_K}=\psi(2g_K) - 4g_K |\mathbb{F}_K|^{g_K}\\&\notag\geq \frac{|\mathbb{F}_K|^{2g_K}}{2g_K} - \left|\frac{1}{2g_K} \sum_{\substack{d|2g_K\\d<2g_K}} \mu\left(\frac{2g_K}{d}\right) |\mathbb{F}_K|^d \right| - 4g_K |\mathbb{F}_K|^{g_K} \\&\geq \frac{|\mathbb{F}_K|^{2g_K}}{2g_K} - \sum_{\substack{d|2g_K\\d<2g_K}} |\mathbb{F}_K|^d  - 4g_K |\mathbb{F}_K|^{g_K}\\&\notag \geq \frac{|\mathbb{F}_K|^{2g_K}}{2g_K} - \sum_{d=1}^{g_K} |\mathbb{F}_K|^d  - 4g_K |\mathbb{F}_K|^{g_K}\\&\notag \geq \frac{|\mathbb{F}_K|^{2g_K}}{2g_K}-(4g_K+2) |\mathbb{F}_K|^{g_K}.\end{align} By (\ref{classnumberlowerbound1}), if $g_K$ is large enough, it holds for any possible value of $|\mathbb{F}_K|$ that \begin{equation} \label{classnumberlowerbound2} h_K \geq  \frac{(|\mathbb{F}_K|-1)|\mathbb{F}_K|^{g_K-1}}{4g_K}. \end{equation} It may be assumed that $g_K > 0$. As $|\mathbb{F}_K| \geq 2$, by application of the logarithm to (\ref{classnumberlowerbound2}), one obtains for $g_K$ large enough that \begin{equation}  \frac{\ln h_K}{g_K \ln |\mathbb{F}_K|} \geq \frac{\ln(|\mathbb{F}_K|-1)}{g_K \ln |\mathbb{F}_K|}+1 - \frac{1}{g_K}-\frac{\ln 4g_K}{g_K \ln |\mathbb{F}_K|} \geq  1 - \frac{1+ \ln 4g_K}{g_K \ln 2}. \end{equation} The result follows. \qed \end{pot1} 

\section{The abelian case}

For this section, let $F$ also denote a congruence function field. Though more involved, the proof of Theorem 2 is also quite basic. The essential steps are as follows. \begin{enumerate} \item Establish the upper bound of Theorem 2 for $K$ with a condition on the growth of the genus via ramification theory and Riemann's inequality; \item Obtain an upper bound for the degree of a finite, geometric, unramified, and abelian extension $H/F$ via global class field theory \cite{Iy}; \item Obtain a lower bound for the degree of the different of a finite abelian extension $K/F$ via higher ramification theory and the Hasse-Arf theorem \cite{Ne1}; \item Derive a contradiction for a sequence that violates the statement of Theorem 2 via the Riemann-Roch theorem, Riemann's hypothesis, and the Riemann-Hurwitz formula. \end{enumerate}

Throughout this section, the notation of Section \ref{section2} is assumed.
The first part of this proof is similar to the method of Madden and Madan in \cite{MaMa}. For a divisor class $C$ of $K$ of degree $n$, one has by Riemann's inequality that $l_K(C) \geq n - g_K + 1$. Thus the number of integral divisors $A_n$ of $K$ of degree $n$ satisfies \begin{equation} \label{degreenbound} A_n \geq h_K\left(\frac{|\mathbb{F}_K|^{n-g_K+1}-1}{|\mathbb{F}_K|-1}\right).\end{equation} Let $s \in \mathbb{R}$ with $s > 1$ and $x \in K \backslash \mathbb{F}_K$. By (\ref{degreenbound}), it follows that \begin{align}\label{Riemanninequality} \zeta_K(s) &= \sum_{n=0}^\infty \frac{A_n}{|\mathbb{F}_K|^{ns}}\geq \sum_{n=g_K}^\infty \frac{A_n}{|\mathbb{F}_K|^{ns}} \geq \sum_{n=g_K}^\infty h_K \left(\frac{|\mathbb{F}_K|^{n-g_K+1}-1}{|\mathbb{F}_K|-1}\right) \frac{1}{|\mathbb{F}_K|^{ns}}\\&\notag= \frac{h_K}{|\mathbb{F}_K|^{{g_K}s}} \sum_{n=g_K}^\infty \frac{|\mathbb{F}_K|^{n-g_K+1}-1}{|\mathbb{F}_K|-1}\frac{1}{|\mathbb{F}_K|^{(n-g_K)s}}=\frac{h_K}{|\mathbb{F}_K|^{{g_K}s}} \zeta_0(s).\end{align} Let $\mathfrak{p} \in \mathbb{P}_0$, and let $\mathfrak{P}_1,...,\mathfrak{P}_r$ be those places of $K$ that lie above $\mathfrak{p}$. For each $i=1,...,r$, let $e(\mathfrak{P}_i|\mathfrak{p})$ denote the ramification index and $f(\mathfrak{P}_i|\mathfrak{p})$ the relative degree of $\mathfrak{P}_i|\mathfrak{p}$. By ramification theory (see for example \cite{Vi}), it holds that \begin{equation}\label{fundamentalidentity} \sum_{i=1}^r e(\mathfrak{P}_i|\mathfrak{p}) f(\mathfrak{P}_i|\mathfrak{p}) = [K:\mathbb{F}_K(x)]. \end{equation} As $K/\mathbb{F}_K(x)$ is a geometric extension, it follows from (\ref{fundamentalidentity}) that \begin{align}\label{euler2} \zeta_K(s) &= \prod_{\mathfrak{P} \in \mathbb{P}_K} \left(1 - \frac{1}{|\mathbb{F}_K|^{d_K(\mathfrak{P})s}}\right)^{-1} \\&\notag\leq \prod_{\mathfrak{p} \in \mathbb{P}_0} \left(1 - \frac{1}{|\mathbb{F}_K|^{d_0(\mathfrak{p})s}}\right)^{-\left[K:\mathbb{F}_K(x)\right]} =\zeta_0(s)^{\left[K:\mathbb{F}_K(x)\right]}. \end{align} By (\ref{Riemanninequality}) and (\ref{euler2}), one obtains that \begin{equation} \label{ratioupperbound} \frac{h_K}{|\mathbb{F}_K|^{{g_K}s}} \leq \zeta_0(s)^{\left[K:\mathbb{F}_K(x)\right]-1}.\end{equation} As $|\mathbb{F}_K| \geq 2$ and $\zeta_0(s) \leq \zeta_{\mathbb{F}_2(T)}(s)$, application of the logarithm to (\ref{ratioupperbound}) yields that \begin{equation}\label{lnzeta} \frac{\ln h_K}{g_K \ln |\mathbb{F}_K|} \leq s + \frac{\left(\left[K:\mathbb{F}_K(x)\right]-1\right)\ln \zeta_{\mathbb{F}_2(T)}(s)}{g_K \ln 2}. \end{equation} Let $\varepsilon$ be fixed and positive, and let $s = 1 + \frac{\varepsilon}{2}$. If the quantity $\left[K:\mathbb{F}_K(x)\right]/g_K$ is chosen to be sufficiently close to zero, it follows from \eqref{lnzeta} that \begin{equation}\notag \frac{\ln h_K}{g_K \ln |\mathbb{F}_K|} < 1 + \varepsilon. \end{equation} The following lemma and the first step of the proof of Theorem 2 have therefore been established.

\begin{lemma} \label{lemma3} Let $x \in K \backslash \mathbb{F}_K$. It holds that \begin{equation} \notag \limsup_{\frac{\left[K:\mathbb{F}_K(x)\right]}{g_K} \rightarrow 0} \: \frac{\ln h_K}{g_K \ln \left|\mathbb{F}_K\right|} \leq 1. \end{equation} \end{lemma} 

By the reciprocity map of global class field theory, a maximal finite, geometric, unramified, and abelian extension of a congruence function field $F$ is of degree $h_F$. From this fact, one obtains the following lemma and the second step of the proof of Theorem 2.

\begin{lemma} \label{lemma4} Let $H/F$ be a finite, geometric, unramified, and abelian extension. It holds that $[H:F] \leq h_F$. \end{lemma}

The third step of the proof of Theorem \ref{theorem2} follows a method known to G. Frey et al. \cite{FrPeSt}. Let $K/F$ be a finite abelian extension. Let $\mathfrak{p} \in \mathbb{P}_F$ and $\mathfrak{P} \in \mathbb{P}_K$ with $\mathfrak{P}|\mathfrak{p}$. For each $n = 0,1,2,...,$ let $G_n(\mathfrak{P}|\mathfrak{p})$ denote the $n$th ramification group of $\mathfrak{P}|\mathfrak{p}$. Also, let $\alpha_{\mathfrak{P}|\mathfrak{p}}$ denote the differential exponent of $\mathfrak{P}|\mathfrak{p}$. First, by ramification theory, one obtains that \begin{equation}\label{Hilbert} \alpha_{\mathfrak{P}|\mathfrak{p}} =\sum_{n=0}^\infty (|G_n(\mathfrak{P}|\mathfrak{p})|-1). \end{equation} Let $k(\mathfrak{P}|\mathfrak{p})$ denote the number of ramification jumps of $\mathfrak{P}|\mathfrak{p}$. By the Hasse-Arf theorem, it follows from (\ref{Hilbert}) that \begin{equation} \label{firstbound} \alpha_{\mathfrak{P}|\mathfrak{p}} \geq \frac{1}{2} k(\mathfrak{P}|\mathfrak{p}) e(\mathfrak{P}|\mathfrak{p}). \end{equation} Second, let $K_\mathfrak{P}$ and $F_\mathfrak{p}$ denote the completion of $K$, respectively $F$, according to $\mathfrak{P}$, respectively $\mathfrak{p}$. Let $\mathfrak{P}$ be identified with its maximal ideal, $\vartheta_\mathfrak{P}$ denote the valuation ring for $\mathfrak{P}$, and $\pi_\mathfrak{P}$ be an element prime for $\mathfrak{P}$. As $K_\mathfrak{P}/F_\mathfrak{p}$ is abelian, the action of $\text{Gal}(K_\mathfrak{P}/F_\mathfrak{p})$ is trivial on each element in the image of each injection \begin{equation}\notag \psi_0 : G_0(\mathfrak{P}|\mathfrak{p})/G_{1}(\mathfrak{P}|\mathfrak{p}) \rightarrow (\vartheta_\mathfrak{P}/\mathfrak{P})^*,\;\;\psi_0(\sigma) = \frac{\sigma(\pi_{\mathfrak{P}})}{\pi_{\mathfrak{P}}}  \end{equation} and, for each $n \in \mathbb{N}$, \begin{equation}\notag \psi_n : G_n(\mathfrak{P}|\mathfrak{p})/G_{n+1}(\mathfrak{P}|\mathfrak{p}) \rightarrow \mathfrak{P}^{n}/\mathfrak{P}^{n+1},\;\;\psi_n(\sigma)= \frac{\sigma(\pi_{\mathfrak{P}})}{\pi_{\mathfrak{P}}} - 1. \end{equation} Identifying $\mathfrak{p}$ with its maximal ideal and denoting by $\vartheta_\mathfrak{p}$ the valuation ring of $\mathfrak{p}$, it follows that \begin{equation}\label{secondbound} e(\mathfrak{P}|\mathfrak{p}) \leq |\vartheta_{\mathfrak{p}}/\mathfrak{p}|^{k(\mathfrak{P}|\mathfrak{p})}. \end{equation} Finally, observing that the fixed field of the product of the inertia groups $G_0(\mathfrak{P}|\mathfrak{p})$ over all $\mathfrak{p} \in \mathbb{P}_F$ is simply the maximal unramified extension of $F$ in $K$, one obtains the following result as a consequence of (\ref{firstbound}) and (\ref{secondbound}).

\begin{lemma} \label{lemma5} Let $K/F$ be a finite abelian extension. Let $H_{K/F}$ denote the maximal unramified extension of $F$ in $K$. It follows that the different $\mathfrak{D}_{K/F}$ satisfies \begin{equation}\notag d_K(\mathfrak{D}_{K/F}) \geq \frac{\left[K:F\right]}{2\ln|\mathbb{F}_F|} \left(\ln \left[K:F\right] - \ln \left[H_{K/F}:F\right] \right). \end{equation} \end{lemma}

\begin{pot2} Consider a sequence $\{K_n\}_{n \in \mathbb{N}}$ with $K_n/F$ a finite abelian extension for each $n \in \mathbb{N}$ and unbounded sequence of genera $\{g_{K_n}\}_{n \in \mathbb{N}}$. Furthermore, suppose that there exists a positive $\delta \in \mathbb{R}$ with, for each $n \in \mathbb{N}$, $\ln h_{K_n}/(g_{K_n} \ln |\mathbb{F}_{K_n}|) \geq 1 + \delta$. Let $x \in F\backslash \mathbb{F}_F$. By Lemma \ref{lemma3}, there exists a positive $\varepsilon \in \mathbb{R}$ with, for each $n \in \mathbb{N}$, \begin{equation} \label{epsilonbound} \frac{\left[K_n:\mathbb{F}_{K_n}(x)\right]}{g_{K_n}}\geq \varepsilon. \end{equation}

Let $s \in \mathbb{C}$, $u = |\mathbb{F}_K|^{-s}$, and $n \in \mathbb{N}$. Let $P_{K_n}(s)$ be defined as in Section \ref{section2}. As noted in (\ref{zetanumerator}), there exist $\omega_1,...,\omega_{2g_{K_n}}$ so that \begin{equation}\label{zetanumerator2} P_{K_n}(s) = \prod_{i=1}^{2g_{K_n}} (1 - \omega_i u). \end{equation} By Riemann's hypothesis, one has for each $i=1,...,2g_{K_n}$ that $|\omega_i| = |\mathbb{F}_{K_n}|^{\frac{1}{2}}$. Also, it is well known \cite{Deu} that $P_{K_n}(0) = h_{K_n}$. From (\ref{zetanumerator2}), one obtains that \begin{equation}\label{classnumberupperbound} h_{K_n} = P_{K_n}(0) = |P_{K_n}(0)| = \prod_{i=1}^{2g_{K_n}} |1 - \omega_i| \leq \left(1 + |\mathbb{F}_{K_n}|^{\frac{1}{2}}\right)^{2g_{K_n}}. \end{equation} It may be assumed for each $n \in \mathbb{N}$ that $g_{K_n} >0$. Application of the logarithm to (\ref{classnumberupperbound}) yields that \begin{equation}\label{finitefieldbound} \frac{\ln h_{K_n}}{g_{K_n}\ln|\mathbb{F}_{K_n}|} \leq \frac{2\ln \left(1 + |\mathbb{F}_{K_n}|^{\frac{1}{2}}\right)}{\ln|\mathbb{F}_{K_n}|}. \end{equation} By (\ref{finitefieldbound}), it follows that the field \begin{equation}\label{E} \mathbb{E} = \prod_{n \in \mathbb{N}} \mathbb{F}_{K_n}\end{equation} is finite. By the definition of $\mathbb{E}$ in ($\ref{E}$), it follows for each $n \in \mathbb{N}$ that the extension $\mathbb{E}K_n/\mathbb{E}F$ is geometric. By the Riemann-Hurwitz formula and Lemmas \ref{lemma4} and \ref{lemma5}, one obtains that \begin{align}\label{riemannhurwitz} \notag \frac{g_{\mathbb{E}K_n}}{\left[\mathbb{E}K_n:\mathbb{E}F\right]}& \geq g_{\mathbb{E}F}-1 + \frac{1}{2\left[\mathbb{E}K_n:\mathbb{E}F\right]} d_{\mathbb{E}K_n}(\mathfrak{D}_{\mathbb{E}K_n/\mathbb{E}F})\\& \geq g_{\mathbb{E}F}-1 + \frac{1}{4 \ln |\mathbb{E}|}(\ln \left[\mathbb{E}K_n:\mathbb{E}F\right] - \ln [H_{\mathbb{E}K_n/\mathbb{E}F}:\mathbb{E}F])\\&\notag \geq g_{\mathbb{E}F}-1 + \frac{1}{4 \ln |\mathbb{E}|}(\ln \left[\mathbb{E}K_n:\mathbb{E}F\right] -  \ln h_{\mathbb{E}F}). \end{align} By basic function field theory, it holds that $\left[\mathbb{E}K_n:\mathbb{E}(x)\right]=\left[K_n:\mathbb{F}_{K_n}(x)\right]$. As the sequence of genera $\{g_{K_n}\}_{n \in \mathbb{N}}$ is unbounded, it follows from (\ref{epsilonbound}) that the sequence $\left\{\left[\mathbb{E}K_n:\mathbb{E}F\right]\right\}_{n \in \mathbb{N}}$ is also unbounded. However, by the Riemann-Roch theorem, one obtains for each $n \in \mathbb{N}$ that $g_{\mathbb{E}K_n} = g_{K_n}$. By (\ref{epsilonbound}) and (\ref{riemannhurwitz}), it follows that the sequence $\left\{\left[\mathbb{E}K_n:\mathbb{E}F\right]\right\}_{n \in \mathbb{N}}$ is bounded. This is a contradiction. The result follows. \qed \end{pot2}

\section{Acknowledgements}

This paper was written while the author was supported by a Vaughn Research Fellowship at Oklahoma State University.





\bibliographystyle{model1b-num-names}
\bibliography{article1references}







\end{document}